\documentclass[12pt,twoside]{article}

\date{}

\usepackage{amssymb}
\usepackage{amsmath}
\usepackage{stmaryrd}

\begin{document}


\centerline{}

\centerline{}

\centerline {\Large{\bf Quasi--bases for Modules over a Commutative Ring}}

\centerline{}


\centerline{\bf {Guang Shi}}

\centerline{}

\centerline{LMIB and School of Mathematics and Systems Science,}

\centerline{Beihang University, Beijing 100191, P.R. China}

\centerline{g\_shi@ss.buaa.edu.cn}








\newtheorem{Theorem}{\quad Theorem}[section]

\newtheorem{Definition}[Theorem]{\quad Definition}

\newtheorem{Corollary}[Theorem]{\quad Corollary}

\newtheorem{Lemma}[Theorem]{\quad Lemma}

\newtheorem{Example}[Theorem]{\quad Example}

\newtheorem{Proposition}[Theorem]{\quad Proposition}

\begin{abstract} In this paper we present the definition of quasi-bases for modules over a ring that is commutative but not necessarily division and discuss properties that guarantee the existence of quasi-bases. Based on this result we further prove that every finitely generated module over $L^{0}(\mathcal{F},K)$ has a quasi-basis, where $K$ is the scalar field of real numbers or complex numbers and $L^{0}(\mathcal{F},K)$ is the algebra of equivalence classes of $K$--valued random variables defined on a probability space $(\Omega,\mathcal{F},P)$.
\end{abstract}

{\bf Mathematics Subject Classification:} 13C99 \\

{\bf Keywords:} quasi--bases, modules, commutative rings, $L^{0}(\mathcal{F},K)$--modules

\section{Introduction}

Since the possible absence of a linear independent subset, it is well known that a unitary module $M$ over a unitary and commutative ring $R$ may not be free if $R$ is not a division ring. Specifically, even $ax=\theta$ ($\theta$ is the null element of $M$) does not imply $a=0$ or $x=\theta$ for $a\in R$ and $x\in M$ because not every non-zero element in $R$ has a multiplicative inverse. To find a group of elements which plays almost the same role as a basis in $M$, one needs to give reasonable generalizations of the classical inverse and linear independence.

The recent progress in random metric theory shed light on finding such generalizations. Random metric theory originated from the theory of probabilistic metric spaces \cite{1, 2, 3, 4, 5, 6, 7, 8}. In the direction of functional analysis, random metric theory has undergone a systematic and deep development \cite{9, 10, 11, 12, 13, 14, 15, 16, 17, 18, 19, 20, 21} since Guo presented the elaborated definition of a random normed module in \cite{10}. Now the theory of random normed modules, random locally convex modules together with their random conjugate spaces, as a random generalization of the theory of classic normed spaces, locally convex spaces and their conjugate spaces, has become a powerful tool for the study of conditional risk measures\cite{16, 22}. Modules over the algebra $L^{0}(\mathcal{F},K)$ of equivalence classes of $K$--valued random variables on some probability space $(\Omega,\mathcal{F},P)$ under the ordinary addition, multiplication and scalar multiplication operations on equivalence classes have been proved to be a kind of proper random generalizations of classic linear spaces over scalar field for random metric theory. Notice that $L^{0}(\mathcal{F},K)$ is not a division ring unless the probability space $(\Omega,\mathcal{F},P)$ is trivial, thus modules over $L^{0}(\mathcal{F},K)$ (briefly, $L^{0}(\mathcal{F},K)$--modules) may not be free. A crucial step in the study of $L^{0}(\mathcal{F},K)$--modules was taken in \cite{14}, where the notion of stratifications were presented. And the algebraic structure of finitely generated $L^{0}(\mathcal{F},K)$--modules was studied in \cite{21}.
In this paper we will extend some results for $L^{0}(\mathcal{F},K)$--modules to those for general modules over a commutative ring.

Generally, let $I:=\{a\in R~|~a^{2}=a\}$ for a unitary and commutative ring $R$,
then $I$ is a lattice if we define $\vee: I\times I \mapsto I$ by $a\vee b=a+b-ab$ and $\wedge:I\times I \mapsto I$ by $a\wedge b=ab$ for any $a,b \in I$. Moreover, if there exists a minimum element $i_{a}$ in $\{i\in I~|~ia=a\}$ under the partial order $\geqslant$ which is induced by $\vee$ and $\wedge$ for each $a\in R$, then a generalized inverse of a non zero element $a$ in $R$ can be defined as the element $a^{-1}$ which satisfies $i_{a^{-1}}=i_{a}$ and $a^{-1}a=i_{a}$. Furthermore, suppose every non zero element in $R$ has a generalized inverse and $M$ is an unitary module over $R$ such that there exists a minimum element $i_{x}$ in $\{i\in I~|~ix=x\}$ for each $x$ in $M$, then it is easy to see that $ax=\theta$ implies $ai_{x}=0$ for any $a\in R$ and $x\in M$. Consequently, a generalized form of linear independence can be defined as: $\{x_{1},x_{2},\cdots,x_{n}\}\subset M$ is said to be linearly independent provided that $\sum_{j=1}^{n}a_{j}x_{j}=\theta$ implies $a_{j}i_{x_{j}}=0,j=1,2,\cdots,n$ for any $a_{1},a_{2},\cdots,a_{n}\in R$. In this paper we will show that if in addition, $R$ has the countable concatenation property and $M$ is finitely generated, i.e. there exists a finite subset $S$ of $M$ such that $M=\mbox{span}S$, then there exists a generalized linearly independent subset $\{x_{1},x_{2},\cdots,x_{n}\}\subset M$ with $i_{x_{1}}\geqslant i_{x_{2}} \geqslant i_{x_{3}} \geqslant\cdots \geqslant i_{x_{n}}$ that can generate $M$, which we called a quasi-basis for $M$. Moreover, we will prove that every finitely generated $L^{0}(\mathcal{F},K)$-module has a quasi-basis.


\section{Main Results}
In this section, $R$ always denotes a commutative ring with identity $e$ and $M$ a unitary module over $R$. Let $I:=\{a\in R~|~a^{2}=a\}$, then $I$ is a lattice if we define $\vee: I\times I \mapsto I$ by $a\vee b=a+b-ab$ and $\wedge:I\times I \mapsto I$ by $a\wedge b=ab$ for any $a,b \in I$. In the sequel of this section we always assume that $I$ is a complete lattice, i.e. every subset $H$ of $I$ has a supremum and an infimum, denoted by $\bigvee H$ and $\bigwedge H$, respectively. In addition, we also assume that there exist subsets $\{a_{n};n\in N\}$ and $\{b_{n};n\in N\}$ of $H$ such that $\bigvee H=\bigvee\{a_{n};n\in N\}$ and $\bigwedge H=\bigwedge\{b_{n};n\in N\}$, where $N$ denotes the set of positive integers. Clearly, if $H$ is directed upwards (downwards), then the above $\{a_{n}\}$ (resp. $\{b_{n}\}$) can be chosen as nondecreasing (correspondingly, nonincreasing). Some of the concepts stated below come from the random metric theory, and here we will restate them in a more general form.

A subset $\{i_{n};n\in N\}$ of $I$ is called a countable partition of $e$ if $\bigvee\{i_{n};n\in N\}=e$ and $i_{m}i_{n}=0$ for any $m,n\in N$ and $m\neq n$. We make the following convention that for any two elements $x$ and $y\in M$, if there exists a countable partition $\{i_{n};n\in N\}$ of $e$ such that $i_{n}x=i_{n}y$ for each $n\in N$, then $x=y$.

\begin{Definition}
A formal sum $\Sigma_{n\in N}i_{n}x_{n}$ is called a countable concatenation of a sequence $\{x_{n}~|~n\in N\}$ in $M$ with respect to a countable partition $\{i_{n};n\in N\}$ of $e$. Moreover, a countable concatenation
$\Sigma_{n\in N}i_{n}x_{n}$ is well defined or $\Sigma_{n\in N}i_{n}x_{n}\in M$
if there is $x\in M$ such that $i_{n}x=i_{n}x_{n}$, $\forall n\in N$.
A subset $G$ of $M$ is called having the countable concatenation property if every countable
concatenation $\Sigma_{n\in N}i_{n}x_{n}$ with $x_{n}\in G$ for each $n\in N$
still belongs to $G$, namely $\Sigma_{n\in N}i_{n}x_{n}$ is well defined and there
exists $x\in G$ such that $x=\Sigma_{n\in N}i_{n}x_{n}$.
\end{Definition}

The following lemma is a restatement of {\cite[Theorem 3.13]{15}}. Here we also give its proof for the convenience of the readers.

\begin{Lemma}[{\cite[Theorem 3.13]{15}}]\label{lem1}
Let $G$ and $H$ be any two nonempty subsets of $M$ such that $G$ and $H$ have the countable concatenation property and $G\cap H=\emptyset$, then there exists an unique element $i_{(G,H)}\in I$ such that the following are satisfied:

\begin{enumerate}
\item $i_{(G,H)}>0$;

\item $iG\cap iH=\emptyset$ for all $i\in I, i\leqslant i_{(G,H)}$ with $i>0$;

\item $iG\cap iH\neq\emptyset$ for all $i\in I, i\wedge i_{(G,H)}=0$ with $i>0$.\\
\end{enumerate}
Here $\leqslant$ is the order induced by $\vee$ and $\wedge$; and $i_{0}>i_{1}$ means $i_{0}\geqslant i_{1}$ and $i_{0}\neq i_{1}$ for $i_{0},i_{1}\in I$.
\end{Lemma}

{\bf Proof}
Let $J=\{i\in I~|~iG\cap iH\neq \emptyset\}$, then for any $i\in J$ we have $i^{\prime}\in J$ for all $i^{\prime}\in I$ with $i^{\prime}\leqslant i$. Moreover, $J$ is directed upwards. Actually, if $i_{1},i_{2}\in J$ then there exist $x_{1}, x_{2}\in G$ and $y_{1},y_{2}\in H$ such that $i_{1}x_{1}=i_{1}y_{1}$ and $i_{2}x_{2}=i_{2}y_{2}$. Let $i_{3}=i_{1}-i_{1}i_{2}$, we have $i_{3}x_{1}+(e-i_{3})x_{2}\in G$, $i_{3}y_{1}+(e-i_{3})y_{2}\in H$ since $G$ and $H$ have the countable concatenation property. Consequently
\begin{eqnarray*}
(i_{1}\vee i_{2})(i_{3}x_{1}+(e-i_{3})x_{2})&=&(i_{1}
-i_{1}i_{2})x_{1}+i_{2}x_{2}\\
&=&(i_{1}
-i_{1}i_{2})i_{1}x_{1}+i_{2}x_{2}\\
&=&(i_{1}
-i_{1}i_{2})i_{1}y_{1}+i_{2}y_{2}\\
&=&(i_{1}
-i_{1}i_{2})y_{1}+i_{2}y_{2}\\
&=&(i_{1}\vee i_{2})
(i_{3}y_{1}+(e-i_{3})y_{2}).\\
\end{eqnarray*}
Hence $i_{1}\vee i_{2}\in J$.

Now let $\{j_{n}~|~n\in N\}$ be a nondecreasing sequence in $J$ such that $\bigvee\{j_{n}~|~n\in N\}=\bigvee J$, then there exist two sequences $\{x_{n}~|~n\in N\}\subset G$ and $\{y_{n}~|~n\in N\}\subset H$ such that $j_{n}x_{n}=j_{n}y_{n}$ for each $n\in N$. Let $j'_{n}=j_{n}-j_{n-1}$ for $n\geqslant 1$, where $j_{0}=0$, then $\Sigma_{n\in N}j'_{n}x_{n}=\Sigma_{n\in N}j'_{n}y_{n}\in G\cap H$ since $G$ and $H$ have the countable concatenation property. Hence $\bigvee J\in J$. Clearly $\bigvee J\neq e$ because $G\cap H=\emptyset$. If we set $i_{(G,H)}=e-\bigvee J$, then it is easy to see that $i_{(G,H)}$ is the desired element.$\Box$

Clearly, a singleton $\{x\}\subset M$ has the countable concatenation property. In the sequel of this section, $i_{x}$ always denotes $i_{(\{x\},\{\theta\})}$ if $x$ is a non-zero element in $M$ and $i_{\theta}=0$, where $\theta$ is the null element of $M$. Specially, $i_{a}:=i_{(\{a\},\{0\})}$ for any $a\in R$ such that $a\neq 0$, and $i_{0}=0$. It is easy to check that $i_{x}$ is the minimum element of $\{i\in I~|~ix=x\}$ for any $x\in M$.

\begin{Definition}
A subset $\{x_{1},x_{2},\cdots,x_{n}\}$ of $M$ is called a quasi-basis if
$M=\mbox{span}\{x_{1},x_{2},\cdots,x_{n}\}$, $i_{x_{1}}\geqslant i_{x_{2}} \geqslant i_{x_{3}} \geqslant\cdots \geqslant i_{x_{n}}$ and for any $\{a_{j}\}_{j=1}^{n}\subset R$, $\Sigma_{j=1}^{n}a_{j}x_{j}=\theta$ implies $a_{j}i_{x_{j}}=0$ for all $j$ such that $1\leqslant j \leqslant n$.
\end{Definition}

Now we can state and prove the main result of this paper as follows:

\begin{Theorem}
If $R$ is a commutative ring with identity such that:
\begin{enumerate}
\item $I=\{a\in R~|~a^{2}=a\}$ is a complete lattice, and for every subset $H$ of $I$ there exist subsets $\{a_{n};n\in N\}$ and $\{b_{n};n\in N\}$ of $H$ such that $\bigvee H=\bigvee\{a_{n};n\in N\}$ and $\bigwedge H=\bigwedge\{b_{n};n\in N\}$;

\item there exists $a^{-1}$ such that $i_{a^{-1}}=i_{a}$ and $a^{-1}a=i_{a}$ for each $a\in R$;

\item $R$ has the countable concatenation property,

\end{enumerate}

then every finitely generated unitary $R$-module $M$ has a quasi-basis; and if $\{x_{1},x_{2},\cdots,x_{n}\}$ and $\{y_{1},y_{2},\cdots,y_{m}\}$ are two quasi-bases for $M$, then $n=m$ and $i_{x_{j}}=i_{y_{j}}$ for $1\leqslant j \leqslant n$.
\end{Theorem}

{\bf Proof}
Suppose $M=\{\Sigma_{j=1}^{l}a_{i}z_{i}~|~a_{i}\in R,1\leqslant j\leqslant l\}$ for some non-zero elements $z_{1},z_{2},\cdots,z_{l}\in M$ and let $z^{(1)}_{1}=z_{1}$. If $a z^{(1)}_{1}=\theta$ for some $a \in R$, then $a i_{x}=0$. In fact, $i_{a}x=a^{-1}(a x)=\theta$ which implies $i_{a}i_{x}=0$. Hence $a i_{x}=a i_{a}i_{x}=0$.

Now suppose for some positive integer $k$ there exists $\{z_{1}^{(k)}, z_{2}^{(k)},\cdots,z_{j}^{(k)}\}$ for some $1\leqslant j\leqslant k$ such that $M_{k}:=\mbox{span}\{z_{1}^{(k)}, z_{2}^{(k)},\cdots,z_{j}^{(k)}\}=\mbox{span}\{z_{1},z_{2},\cdots,z_{l_{k}}\}$ for some $l_{k}\geqslant k$, $i_{z_{1}^{(k)}}\geqslant i_{z_{2}^{(k)}} \geqslant i_{z_{2}^{(k)}} \geqslant\cdots \geqslant i_{z_{j}^{(k)}}$ and for any $\{a_{p}\}_{p=1}^{j}\subset R$, $\Sigma_{p=1}^{j}a_{p}z_{p}^{(k)}=\theta$ implies $a_{p}i_{z_{p}^{(k)}}=0$ for all $p$ such that $1\leqslant p \leqslant j$. If $M_{k}=M$, then $\{z_{1}^{(k)}, z_{2}^{(k)},\cdots,z_{j}^{(k)}\}$ is a quasi-basis for $M$; else let $l_{k+1}$ be the minimum integer such that $z_{l_{k+1}}\notin M_{k}$. Clearly $M_{k}$ has the countable concatenation property thus $i_{(\{z_{l_{k+1}}\},M_{k})}>0$, and if $\Sigma_{p=1}^{j}b_{p}z_{p}^{(k)}+b i_{(\{z_{l_{k+1}}\},M_{k})}z_{l_{k+1}}=\theta$ for some $b_{1},b_{2},\cdots,b_{p},b\in R$, then $$i_{b} i_{(\{z_{i_{k+1}}\},M_{k})}z_{l_{k+1}}=b^{-1}b i_{(\{z_{l_{k+1}}\},M_{k})}z_{l_{k+1}}=-b^{-1}\Sigma_{p=1}^{j}b_{p}z_{p}^{(k)}\in M_{k},$$
By the chosen of $i_{(\{z_{i_{k+1}}\},M_{k})}$ we have $i_{b} i_{(\{z_{i_{k+1}}\},M_{k})}=0$, which implies
$b i_{(\{z_{i_{k+1}}\},M_{k})}=0$ and consequently implies $b_{p}i_{z_{p}^{(k)}}=0$ for $1\leqslant p \leqslant j$.

Moreover, it is easy to check that $\mbox{span}\{z_{1}^{(k)}, z_{2}^{(k)},\cdots,z_{j}^{(k)}, i_{(\{z_{l_{k+1}}\},M_{k})}z_{l_{k+1}}\}=\mbox{span}\{z_{1},z_{2},\cdots,z_{l_{k+1}}\}$. Let $i_{1}=i_{(\{z_{l_{k+1}}\},M_{k})}i_{z_{1}^{(k)}}$ and $z_{1}^{(k+1)}=z_{1}^{(k)}+(i_{(\{z_{l_{k+1}}\},M_{k})}-i_{1})z_{l_{k+1}}$. If $i_{1}=0$, let $z_{p}^{(k+1)}=z_{p}^{(k)}$ for $2\leqslant p\leqslant j$; else let $i_{2}=i_{1}i_{z_{2}^{(k)}}$ and
$z_{2}^{(k+1)}=z_{2}^{(k)}+(i_{1}-i_{2})z_{l_{k+1}}$. Likewise, if $i_{q}=0$ for some $2\leqslant q\leqslant j-1$, let $z_{p}^{(k+1)}=z_{p}^{(k)}$ for $q< p\leqslant j$; else let $i_{q+1}=i_{q}i_{z_{q+1}^{(k)}}$ and
$z_{q+1}^{(k+1)}=z_{q+1}^{(k)}+(i_{q}-i_{q+1})z_{l_{k+1}}$. Finally, if $i_{j}\neq 0$ let $z_{j+1}^{(k+1)}=i_{j}z_{l_{k+1}}$. it is easy to verify that $\mbox{span}\{z_{1}^{(k+1)}, z_{2}^{(k+1)},\cdots,z_{r}^{(k+1)}\}=\mbox{span}\{z_{1},z_{2},\cdots,z_{l_{k+1}}\}$, $i_{z_{1}^{(k+1)}}\geqslant i_{z_{2}^{(k+1)}} \geqslant i_{z_{3}^{(k+1)}} \geqslant\cdots \geqslant i_{z_{r}^{(k+1)}}$ and for any $\{c_{p}\}_{p=1}^{r}\subset R$, $\Sigma_{p=1}^{r}c_{p}z_{p}^{(k)}=\theta$ implies $c_{p}i_{z_{p}^{(k)}}=0$ for all $p$ such that $1\leqslant p \leqslant r$; where $r=j+1$ if $i_{j}\neq 0$ or else $r=j$.

By induction we could obtain a quasi-basis for $M$. Now we turn to prove the later part of the theorem. Let $i_{0}=i-i_{x_{1}}$, $i_{j}=i_{x_{j}}-i_{x_{j+1}}$ for $1\leqslant j< n$, $i_{n}=i_{x_{n}}$, $i'_{0}=i-i_{y_{1}}$, $i'_{j}=i_{y_{j}}-i_{y_{j+1}}$ for $1\leqslant l< m$ and $i'_{m}=i_{y_{m}}$; then it is easy to check that $i_{0}M:=\{i_{0}x~|~x\in M\}=\{\theta\}$ and $\{i_{j}x_{1},i_{j}x_{2},\cdots,i_{j}x_{j}\}$ is a basis for the module $i_{j}M$ over the ring $i_{j}R:=\{i_{j}a~|~a\in R\}$ for each $j$ satisfying $1\leqslant j\leqslant n$ and $i_{j}>0$. Thus $i_{j}M$ is a free module of rank $j$ over the ring $i_{j}R$ for each $j$  satisfying $0\leqslant j\leqslant n$ and $i_{j}>0$. Likewise $i'_{l}M$ is a free module of rank $l$ over the ring $i'_{l}R$ for each $l$  satisfying $0\leqslant l\leqslant m$ and $i'_{l}>0$.

Suppose $i_{j}i'_{l}>0$ for some $0\leqslant j\leqslant n$ and $0\leqslant l\leqslant m$, then it is easy to check that $i_{j}i'_{l}M$ is a free module of rank $j$ and also $l$ over the ring $i_{j}i'_{l}R$. Since $i_{j}i'_{l}R$
is a commutative ring with identity $i_{j}i'_{j}$, it follows from \cite[Chapter 4, Corollary 2.12]{24} that $i_{j}i'_{l}E$ has the invariant dimension property, which implies $j=l$. Notice that $$i_{j}i'_{j}=i_{j}\sum_{l=0}^{m}i'_{l}=i_{j}e=i_{j},$$
thus $i_{j}<i'_{j}$ for $0\leqslant j \leqslant n$. Similarly, we have $i'_{l}<i_{l}$ for $0\leqslant l \leqslant m$, hence $m=n$ and $i_{x_{j}}=i_{y_{j}}$ for $1\leqslant j \leqslant n$, which completes our proof.$\Box$

Now we turn to prove that every finitely generated $L^{0}(\mathcal{F},K)$-module has a quasi-basis.
In the sequel of this paper, $\tilde{I}_{A}$ denotes the equivalence class in $L^{0}(\mathcal{F},K)$ determined by the characteristic function of an $\mathcal{F}$-measurable subset $A$ of $\Omega$; and $\mathcal{I}=\{\tilde{I}_{A}~|~A\in \mathcal{F}\}$. It is easy to check that $\mathcal{I}=\{\xi\in  L^{0}(\mathcal{F},K)~|~\xi^{2}=\xi\}$, $\tilde{I}_{A}\vee\tilde{I}_{B}:=\tilde{I}_{A\cup B}$ and $\tilde{I}_{A}\wedge\tilde{I}_{B}:=\tilde{I}_{A\cap B}$ for any $A,B\in \mathcal{F}$, thus $\mathcal{I}$ is a lattice under $\vee$ and $\wedge$. The following proposition is easy to see from \cite{23}.

\begin{Proposition}
$\mathcal{I}$ is a complete lattice, i.e. every subset $\mathcal{H}$ of $\mathcal{I}$ has a supremum and an infimum, denoted by $\bigvee \mathcal{H}$ and $\bigwedge \mathcal{H}$, respectively; and there exist subsets $\{a_{n};n\in N\}$ and $\{b_{n};n\in N\}$ of $\mathcal{H}$ such that $\bigvee \mathcal{H}=\bigvee\{a_{n};n\in N\}$ and $\bigwedge \mathcal{H}=\bigwedge\{b_{n};n\in N\}$, where $N$ stands for the set of positive integers. Furthermore if $\mathcal{H}$ is directed upwards (downwards) then $\{a_{n};n\in N\}$ can be chosen as nondecreasing (resp. $\{b_{n};n\in N\}$ can be chosen as nonincreasing).
\end{Proposition}

Every countable partition of $\tilde{I}_{\Omega}$ can be represented as $\{\tilde{I}_{A_{n}}; n\in N\}$, where $\{A_{n}; n\in N\}$ is a countable partition of $\Omega$ to $\mathcal{F}$. Clearly, $L^{0}(\mathcal{F},K)$ has the countable concatenation property. Moreover, for any $\xi\in L^{0}(\mathcal{F},K)$, $\xi^{-1}$ is exactly the equivalence class of the $\mathcal{F}$-measurable function $(\xi^{0})^{-1}:\Omega\rightarrow K$ defined by
$$
(\xi^{0})^{-1}(\omega)=\left \{
\begin{array}{ll}
(\xi^{0}(\omega))^{-1}, &\mbox{if~}\xi^{0}(\omega)\neq 0;\\
0, &\mbox{otherwise},
\end{array}
\right.
$$ where $\xi^{0}$ is an arbitrarily chosen representative of $\xi$. It is clear that
and $\xi \cdot \xi^{-1}=\tilde{I}_{\{\omega\in \Omega~|~\xi^{0}(\omega)\neq 0\}}=i_{\xi}$. Thus we have the following theorem:

\begin{Theorem}
If $E$ is a finitely generated $L^{0}(\mathcal{F},K)$-module, then $E$ has a quasi-basis. Moreover, if $\{x_{1},x_{2},\cdots,x_{n}\}$ and $\{y_{1},y_{2},\cdots,y_{m}\}$ are two quasi-basis for $E$, then $m=n$ and $i_{x_{i}}=i_{y_{i}}$ for $1\leqslant i \leqslant n$.
\end{Theorem}




\begin{thebibliography}{99}


\bibitem{1} B. Schweizer and A. Sklar, Probabilistic Metric
Spaces, Elsevier, New York, 1983.

\bibitem{2} B. Schweizer and A. Sklar, Probabilistic Metric Spaces,
Dover Publications, New York, 2005.

\bibitem{3} C. Alsina, B. Schweizer, A. Sklar, On the definition of a probabilistic normed space, Aequationes Math. 46 (1993) 91-98.

\bibitem{4} C. Alsina, B. Schweizer, C. Sempi, A. Sklar, On the definition of a probabilistic inner product space, Rend. Mat. 17 (1997) 115-127.

\bibitem{5} C. Alsina, B. Schweizer, A. Sklar, Continuity property of probabilistic norms, J. Math. Anal. Appl. 208 (1997) 446-452.

\bibitem{6} B. Lafuerza-Guill\'{e}n, J.A. Rodr\'{i}guez-Lallena, C. Sempi, A study of boundedness in probabilistic normed spaces, J. Math. Anal. Appl. 232 (1999) 183-196.

\bibitem{7} B. Lafuerza-Guill\'{e}n, C. Sempi, Probabilistic norms and convergence of random variables, J. Math. Anal. Appl. 280 (2003) 9-16.

\bibitem{8} C. Sempi, A short and partial history of probabilistic normed spaces, Mediterr. J. Math. 3 (2006) 283-300.


\bibitem{9} T.X. Guo, Extension theorems of continuous random linear
operators on random domains, J. Math. Anal. Appl. 193(1) (1995) 15-27.

\bibitem{10} T.X. Guo, Some basic theories of random normed linear
spaces and random inner product spaces, Acta Anal. Funct.
 Appl. 1(2) (1999) 160-184.

\bibitem{11} T.X. Guo, The Radon-Nikod\'ym property of conjugate spaces and the w$\ast$-equivalence theorem for w$\ast$-measurable functions, Sci. China Ser. A 39 (1996) 1034-1041.

\bibitem{12} T.X. Guo, Representation theorems of the dual of
Lebesgue-Bochner function spaces, Sci. China Ser. A 43 (2000)
234-243.

\bibitem{13} T.X. Guo, Several applications of the theory of random
conjugate spaces to measurability problems, Sci. China Ser. A 50 (2007) 737-747.

\bibitem{14} T.X. Guo, The relation of Banach-Alaoglu theorem and
Banach-Bourbaki-Kakutani-$\check{\textmd{S}}$mulian theorem in
complete random normed modules to stratification structure, Sci.
China Ser. A 51 (2008) 1651-1663.

\bibitem{15} T.X. Guo, Relations between some basic results derived from two
kinds of topologies for a random locally convex module, J. Funct. Anal. 258 (2010)
3024-3047.

\bibitem{16} T.X. Guo, Recent progress in random metric theory and its applications
to conditional risk measures, Sci. China Ser. A 71 (2009) 3794-3804.

\bibitem{17} T.X. Guo and S.B. Li, The James theorem in complete
random normed modules, J. Math. Anal. Appl. 308 (2005) 257-265.

\bibitem{18} T.X. Guo and S.L. Peng, A characterization for an $L(\mu,K)$-topological
module to admit enough canonical module homomorphisms, J. Math. Anal.
Appl. 263 (2001) 580-599.

\bibitem{19} T.X. Guo, H.X. Xiao and X.X. Chen, A basic strict separation
theorem in random locally convex modules, Nonlinear Anal. 71 (2009)
3794-3804.

\bibitem{20} T.X. Guo and Z.Y. You, The Riesz's representation theorem in complete random inner product modules and its applications, Chinese Ann. Math. Ser. A 17 (1996) 361-364.

\bibitem{21} T.X. Guo and G. Shi, The algebraic structure of finitely generated $L^{0}(\mathcal{F},K)$-modules and Helly theorem in random normed modules, J. Math. Anal. Appl. 381 (2011) 833-842.


\bibitem{22} D. Filipovi\'{c}, M. Kupper and N. Vogelpoth, Separation
and duality in locally $L^{0}$-convex modules, J. Funct. Anal. 256 (2009)
3996--4029.


\bibitem{23} N. Dunford and J. T. Schwartz, Linear
Operators(\uppercase \expandafter {\romannumeral 1}), Interscience,
New York, 1957.


\bibitem{24} T.W. Hungerford, Algebra, Springer, New York, 1974.

\end{thebibliography}
\end{document}